\input amstex \documentstyle {amsppt}
\magnification=\magstep1
\define\cA{{\Cal A}}

\define\cB{{\Cal B}}

\define\cF{{\Cal F}}
\define\cT{{\Cal T}}
\define\ee{\text{E}}
\define\RR{\Bbb R}
\define\NN{\Bbb N}
\define\Nn{\NN_0}

\define\CC{\Bbb C}
\define\HH{{\Cal H}}

\define\ff{\varphi}

\define\lb{\lbrack}
\define\rb{\rbrack}
\define\tr{\operatorname{tr}}
\define\intdotsint{\operatornamewithlimits{\idotsint}}
\define\wlim{\operatornamewithlimits{w-lim}}

\define\id{\text{{\bf 1}}}

\define\la{\langle}
\define\ra{\rangle}

\define\symb{\Psi}
\define\sq{\Gamma}
\define\sqq{\sq_q}
\define\HR{\Cal H}
\define\HC{\HR_{\CC}}
\define\oo{\omega}
\define\cFq{\cF_q}
\define\sqH{\sqq(\HR)}
\define\fen{f_1\otimes\dots\otimes f_n}
\define\vN{\text{vN}}
\define\lmu{\lim_{m\to\infty}}
\define\lku{\lim_{k\to\infty}}
\define\odo{\otimes\dots\otimes}
\define\sqT{\sqq(T)}

\define\cFH{\cF_q(\HH)}
\define\cFT{\cF_q(T)}
\define\LzH{L^2_q(\HH)}

\define\LuH{L_q^\infty(\HH)}

\define\cFf{\cF^{finite}(\HH)}
\define\cFfs{\cF^{finite}(\HH')}
\define\cFHs{\cF_q(\HH')}
\define\cK{\Cal K}
\define\cKC{\cK_{\CC}}

\define\Xt{(X_t)_{t\in T}}
\define\Xtk{(\tilde X_t)_{t\in T}}
\define\XBM{(X_t^{qBM})_{t\in \lb 0,\infty)}}

\define\XBB{(X_t^{qBB})_{t\in \lb 0,1\rb}}
\define\XOU{(X_t^{qOU})_{t\in \RR}}

\define\Aff{(\cA,\ff)}
\define\cAbt{\cA_{t\rb}}
\define\cAgt{\cA_{\lb t\rb}}
\define\cAbs{\cA_{s\rb}}
\define\cAgs{\cA_{\lb s\rb}}
\define\cHbt{\HR_{t\rb}}
\define\cHgt{\HR_{\lb t\rb}}
\define\cHgs{\HR_{\lb s\rb}}
\define\Span{\text{span}}
\define\cPbt{P_{t\rb}}
\define\cPbs{P_{s\rb}}
\define\Phq{\Phi_q}
\define\PhqH{\Phq(\HR)}
\define\PhqT{\Phq(T)}
\define\Hnq{H_n^{(q)}}
\define\spec{\text{spect}}
\define\Kst{\cK_{s,t}}
\define\Kstq{\cK^{(q)}_{s,t}}
\define\ck{k}
\define\kstq{\ck^{(q)}_{s,t}}
\define\kst{\ck_{s,t}}
\define\kstn{\ck^{(0)}_{s,t}}
\define\lst{{\lambda_{s,t}}}
\define\lstn{{\lambda_{s,t}^n}}
\define\lstq{{\lambda_{s,t}^2}}
\define\ls{{\lambda_{s}}}
\define\lt{{\lambda_{t}}}
\define\ltn{{\lambda_{t}^n}}
\define\ltq{{\lambda_{t}^2}}
\define\lsq{{\lambda_{s}^2}}
\define\ftn{f^{\otimes n}}
\define\cKqt{\cK_t^{(q)}}
\define\ckqt{\ck_t^{(q)}}
\topmatter
\title
$q$-Gaussian processes:\\
non-commutative and classical aspects
\endtitle
\rightheadtext{$q$-GAUSSIAN PROCESSES}
\author
Marek Bo\D zejko$^*$, Burkhard K\"ummerer, Roland Speicher$^{**}$
\endauthor
\leftheadtext{M. BO\D ZEJKO, B. K\"UMMERER, R. SPEICHER}
\address
Instytut Matematyczny,
Uniwersytet Wroc\l awski, Plac Grunwaldzki 2/4,
50-384 Wroc\l aw, Poland 
\endaddress
\email
bozejko\@math.uni.wroc.pl
\endemail
\address
Mathematisches Institut A,
Pfaffenwaldring 57, D-70569 Stuttgart,
Germany
\endaddress
\email
kuem\@mathematik.uni-stuttgart.de
\endemail
\address
Institut f\"ur Angewandte Mathematik, Universit\"at Heidelberg,
Im Neuenheimer Feld 294, D-69120 Heidelberg, Germany
\endaddress
\email
roland.speicher\@urz.uni-heidelberg.de
\endemail
\abstract
We examine, for $-1<q<1$, 
 $q$-Gaussian processes, i.e. families of operators 
(non-commutative random variables) $X_t=a_t+a_t^*$ -- where
the $a_t$ fulfill the $q$-commutation relations 
$a_sa_t^*-qa_t^*a_s=c(s,t)\cdot \id$ for some covariance
function $c(\cdot,\cdot)$ -- equipped with
the vacuum expectation state. We show that there is a $q$-analogue
of the Gaussian functor 
of second quantization behind these processes and that this 
structure can be used to translate questions on $q$-Gaussian
processes into corresponding (and much simpler) questions in
the underlying Hilbert space. In particular, we use this idea to
show that a large class of $q$-Gaussian processes possess a
non-commutative kind of Markov property, which ensures that there
exist classical versions of these non-commutative
processes. This answers an old
question of Frisch and Bourret \cite{FB}.
\endabstract
\thanks
$^*$ Partially supported by Polish National Grant, KBN 4233
\endthanks
\thanks
$^{**}$ Supported by a Heisenberg fellowship from the DFG
\endthanks
\thanks
We thank Philippe Biane for stimulating discussions and
remarks.
\endthanks
\endtopmatter
\document
\heading
{\bf Introduction}
\endheading
What we are going to call $q$-Gaussian processes was essentially
introduced in a remarkable paper by Frisch and Bourret
\cite{FB}. Namely, they considered generalized commutation
relations given by operators $A(t)$ and a vacuum vector $\Psi_0$
with
$$A(t)A^*(t')-q A^*(t')A(t)=\Gamma(t,t')\id$$
and
$$A(t)\Psi_0=0$$
for some real covariance function 
$\Gamma$ (i.e. positive definite
function). The aim of the authors was to study the probabilistic
properties of the \lq parastochastic' process
$M(t)=A(t)+A^*(t)$.
\par
The basic problems arising in this context were the following
two types of questions:
\roster
\item"(I)" (realization problem) \newline
Do there exist operators on some Hilbert space and a corresponding
vacuum vector in this Hilbert space 
which fullfill the above relations,
i.e. are there non-commutative realizations of the $q$-Gaussian 
processes.
\item"(II)" (random representation problem) \newline
Are these non-commutative processes of a classical relevance, i.e.
do there exist classical versions of the $q$-Gaussian processes (in the
sense of coinciding time-ordered correlations, see our Def. 4.1.)
\endroster
Frisch and Bourret could give the following partial answers to these
questions.
\roster
\item"(I)"
For $q=\pm 1$ the realization is of course given by the
Fock space realization of the bosonic/fermionic relations.
The case $q=0$ was realized by creation and annihilation operators
on the full Fock space (note that this was before the introduction
of the Cuntz algebras and their extensions \cite{Cun,Eva}).
For other values of $q$ the realization problem remained open.
\item"(II)" The $q=1$ processes are nothing but the Fock space
representations of the classical Gaussian processes. For $q=-1$
a classical realization by a dichotomic Markov process could be
given for the special case of exponential covariance 
$\Gamma(t,t')=\exp(-\vert t-t'\vert)$.
A classical realization for $q=0$ could not be found, but the authors
were able to show that there is an interesting representation
in terms of Gaussian random matrices.
\endroster
The authors started also the investigation of parastochastic 
equations (i.e. 
the coupling of parastochastic processes to other systems), but
 -- probably because of the open problem on the mere existence
and classical relevance of these $q$-processes -- there was
apparently no further work in this direction and the paper of Frisch
and Bourret felt into oblivion.
\par
Starting with \cite{AFL} there has been another and independent 
approach to non-commutative probability theory. This wide and
quite inhomogenous field 
-- let us just mention as two
highlights the quantum stochastic
calculus of Hudson-Parthasarathy
\cite{HP} and the free probability theory
of Voiculescu \cite{VDN} -- is now known under the name of 
\lq quantum probability'. 
At least some of the fundamental motivation for undertaking
such investigations can be compared with the two basic questions
of Frisch and Bourret:
\roster
\item"(I)" Non-commutative probability theory is meant as a 
generalization of classical probability theory to the 
description of quantum systems. Thus first of all their objects are
operators on some Hilbert spaces having a meaning as non-commutative
anlogues of the probabilistic notions of random variables, stochastic
processes, etc.
\item"(II)" In many
investigations in this area one also 
tries to 
establish connections between non-commutative and classical
concepts. The aim of this is twofold. On one side, one hopes
to get a better understanding of classical problems by
embedding them into a bigger non-commutative context.
Thus, e.g., the Az\'ema martingale, although classically not
distinguished within the class of all martingales, behaves in some
respects like a Brownian motion
\cite{Par1}. The non-commutative \lq explanation' for
this fact comes from the observation of Sch\"urmann \cite{Sch}
that this martingale is one component of a non-commutative
process with independent increments. 
In the other direction, one hopes to get a classical picture
(featuring trajectories)
of some aspects of quantum problems. Of course, a total reduction
to classical concepts is in general not possible, but partial aspects
may sometimes allow a classical interpretation.
\endroster
It was in this context of quantum probability where two of the
present authors \cite{BSp1} reintroduced the $q$-relations --
without knowing of, but much in the same spirit as \cite{FB}.
Around the same time the $q$-relations were also proposed by
Greenberg \cite{Gre} as an example for particles with
\lq infinite statistics'. 
\par
The main progress in connection with this renewed interest was the
solution of the realization problem
of Frisch and Bourret. There exist now
different proofs for the existence of the Fock representation
of the $q$-relations for all $q$ with $-1\leq q\leq 1$
\cite{BSp1,Zag,Fiv,Spe1,BSp3,YW}.
\par
In \cite{NSp}, the idea of Frisch and Bourret to use the $q$-relations
as a model for a generalized noise was pursued further and the 
Greens function for such dynamical problems could be calculated
for one special choice of the covariance function -- namely for
the case of the exponential
covariance. We will call this special $q$-process in the following
$q$-Ornstein-Uhlenbeck process.
It soon became clear that the special status of the exponential
covariance is 
connected with some kind of (non-commutative) Markovianity -- as we
will see the $q$-Ornstein-Uhlenbeck process is the only stationary
$q$-Gaussian Markov process. But using the general theory of
K\"ummerer on non-commutative stationary
Markov processes \cite{K\"um1,K\"um2}
this readily implies the existence of a classical
version (being 
itself a classical Markov process) of the $q$-Ornstein-Uhlenbeck
process. Thus we got a positive solution of the random representation
problem of Frisch and Bourret in this case. However, 
the status of the other
$q$-Gaussian processes, in particular $q$-Brownian motion, remained
unclear.
\par
Motivated by our preliminary results, Biane \cite{Bia1} (see
also \cite{Bia2,Bia3}) undertook
a deep and beautiful analysis of the free ($q=0$) case and showed
the remarkable result that all processes with free increments are
Markovian and thus possess classical versions (with a
quite explicit calculation rule for the corresponding transition
probabilities). This includes in particular the case of free
Brownian motion.
\par
Inspired by this work we could extend our investigations from the
case of the $q$-Ornstein-Uhlenbeck process to all $q$-Gaussian
processes. The results are presented in this paper.
\par 
Up to now there is only one strategy for establishing the existence of
a classical version of a non-commutative process, namely by
showing that the process is Markovian. That this implies
the existence of a classical version follows by general arguments,
the main point is to show that we have this property in the concrete case.
Whereas Biane could use the quite developed theory of freeness \cite{VDN} to
prove Markovianity for processes with free increments, there
is at the moment (and probably also in
the future \cite{Spe2}) no kind of
$q$-freeness for general $q$. Thus another feature of our
considered class of processes is needed to attack the problem
of Markovianity. It is the aim of this paper to
convince the reader of the fact that 
the $q$-analogue of Gaussianity will do this job.
\par
The essential idea of Gaussianity is that one can pull back all 
considerations 
from the measure theoretic (or, in the non-commutative frame,
from the operator algebraic) level to an underlying Hilbert space, thus
in the end 
one essentially has to deal with linear problems. The main point is
that this transcription between the linear and the
algebraic level exists in a consistent way. The best way to
see and describe this is by presenting a functor 
(\lq second quantization') which translates
the Hilbert space properties into operator algebraic properties.
Our basic considerations will therefore be on the existence and 
nice properties of the $q$-analogue of this functor. Having this
functor, the rest is mainly linear theory on Hilbert space level. 
It turns out that all relevant questions on our $q$-Gaussian processes
can be characterized totally in terms of the corresponding
covariance function. In particular, it becomes quite easy to
decide whether such a process is Markovian or not.
\par
The paper is organized as follows.
In Sect. 1 we remind of some basic facts about the $q$-Fock
space and its relevant operators. Furthermore we collect 
in this section the
needed combinatorial results, in particular on $q$-Hermite
polynomials.
Sect. 2 is devoted to the presentation of the functor 
$\Gamma_q$ of
second quantization. The main results (apart from the existence
of this object) are the facts that the associated von Neumann
algebras are in the infinite dimensional case 
non-injective $II_1$-factors and
that the functor maps contractions into completely positive maps.
Having this $q$-Gaussian functor the definition and investigation
of properties of $q$-Gaussian processes (like Markovianity or martingale
property)
is quite canonical and parallels the classical case. Thus our
presentation of these aspects, in Sect. 3, will be quite condensed.
Sect. 4 contains the classical interpretation of the $q$-Gaussian
Markov processes. As pointed out 
above general arguments ensure
the existence of classical versions for these
processes. But we will see that we can also derive quite concrete
formulas for the corresponding
transition probabilities.
\heading
{\bf 1. Preliminaries on the $q$-Fock space}
\endheading
Let $q\in(-1,1)$ be fixed in the following.
\par
For a complex Hilbert space $\HH$ we define its $q$-Fock space
$\cFH$ as follows: Let $\cFf$ be the linar span of vectors of the
form $\fen\in\HH^{\otimes n}$ 
(with varying $n\in\Nn$), where we put
$\HH^{\otimes 0}\cong\CC\Omega$ for some distinguished vector
$\Omega$, called vacuum. On $\cFf$ we consider the sesquilinear
form $\la\cdot,\cdot\ra_q$ given by sesquilinear extension of
$$\la\fen,g_1\odo g_m\ra_q:=
\delta_{nm}\sum_{\pi\in S_n}q^{i(\pi)}\la f_1,g_{\pi(1)}\ra\dots
\la f_n,g_{\pi(n)}\ra,$$
where $S_n$ denotes the symmetric group of permutations 
of $n$ elements and $i(\pi)$ is the number of inversions of the
permutation $\pi\in S_n$ defined by
$$i(\pi):=\#\{(i,j)\mid 1\leq i<j\leq n, \pi(i)>\pi(j)\}.$$
Another way to describe $\la\cdot,\cdot\ra_q$ is by introducing
the operator $P_q$ on $\cFf$ by linear extension of
$$\align
P_q\Omega&=\Omega\\
P_q\fen&=\sum_{\pi\in S_n}q^{i(\pi)}f_{\pi(1)}\odo f_{\pi(n)}.
\endalign$$
Then we can write
$$\la\xi,\eta\ra_q=\la\xi,P_q\eta\ra_0\qquad (\xi,\eta\in\cFf),$$
where $\la\cdot,\cdot\ra_0$ is the scalar product on the usual
full Fock space 
$$\cF_0(\HH)=\bigoplus_{n\geq 0}\HH^{\otimes n}.$$
One of the main results of \cite{BSp1} (see also \cite{BSp3,Fiv,Spe1,Zag})
was the strict 
positivity of $P_q$, i.e. $\la\xi,\xi\ra_q>0$ for 
$0\not=\xi\in\cFf$. 
This allows the following
definitions.
\definition{1.1. Definitions}
1) The {\it $q$-Fock space} $\cFH$ is the completion of $\cFf$
with respect to $\la\cdot,\cdot\ra_q$.\newline
2) Given $f\in\HH$, we define the {\it creation operator}
$a^*(f)$ and the {\it annihilation operator} $a(f)$ on $\cFH$ by
$$\align
a^*(f)\Omega&=f\\
a^*(f)\fen&=f\otimes\fen\endalign$$
and
$$\align
a(f)\Omega&=0\\
a(f)\fen&=\sum_{i=1}^nq^{i-1}
\la f,f_i\ra f_1\odo\check f_i\odo f_n
,\endalign$$
where the symbol $\check f_i$ means that $f_i$ has to be
deleted in the tensor.
\enddefinition
\remark{1.2. Remark}
The operators $a(f)$ and $a^*(f)$ are bounded operators on
$\cFH$ with
$$\Vert a(f)\Vert_q=\Vert a^*(f)\Vert_q=
\cases {\Vert f\Vert}/{\sqrt{1-q}},& 0\leq q<1\\
\Vert f\Vert, & -1< q\leq 0,\endcases$$
and they are adjoints of each other with respect to our
scalar product $\la\cdot,\cdot\ra_q$.
Furthermore, they fulfill the $q$-relations
$$a(f)a^*(g)-qa^*(g)a(f)=\la f,g\ra\cdot\id\qquad (f,g\in\HH).$$ 
\endremark
\remark{1.3. Notation}
For a linear operator $T:\HH\to\HH'$ between two complex
Hilbert spaces we denote by
$\cF(T):\cFf\to\cFfs$ 
the linear extension of
$$\align
\cF(T)\Omega&=\Omega\\
\cF(T)\fen&=(Tf_1)\odo (Tf_n).\endalign$$
In order to keep the notation simple we denote the vacuum for
$\HH$ and the vacuum for $\HH'$ by the same symbol $\Omega$.
\endremark
It is clear that $\cF(T)$ can be extended to a bounded operator
$\cF_0(T):\cF_0(\HH)\to\cF_0(\HH')$
exactly if $T$ is a contraction, i.e. if $\Vert T\Vert\leq 1$.
The following lemma ensures that the same is true for all
other $q\in(-1,1)$, too.
\proclaim{1.4. Lemma}
Let $\cT:\cFf\to\cFfs$ be a linear operator which fulfills
$P_q'\cT=\cT P_q$, where $P_q$ and $P_q'$ are the operators
on $\cFf$ and $\cFfs$, respectively, which define the respective
scalar product $\la\cdot,\cdot\ra_q$. Then one has
$\Vert\cT\Vert_q=\Vert\cT\Vert_0$. Hence, if $\Vert\cT\Vert_0<
\infty$ then $\cT$ can, for each $q\in(-1,1)$, be extended 
to a bounded operator from
$\cFH$ to $\cFHs$.
\endproclaim
\demo{Proof}
Let $\xi\in\cFf$. Then
$$\align
\Vert \cT\xi\Vert_q^2&=\la\cT\xi,\cT\xi\ra_q\\
&=\la\cT\xi,P_q'\cT\xi\ra_0\\
&=\la P_q^{1/2}\xi,\cT^*\cT P_q^{1/2}\xi\ra_0\\
&\leq \Vert\cT^*\cT\Vert_0\thinspace \la P_q^{1/2}\xi,P_q^{1/2}\xi
\ra_0\\
&=\Vert\cT^*\cT\Vert_0\thinspace \Vert\xi\Vert_q^2,
\endalign$$
which implies
$$\Vert\cT\Vert^2_q\leq\Vert\cT^*\cT\Vert_0\leq\Vert\cT^*\Vert_0
\thinspace\Vert\cT\Vert_0=\Vert\cT\Vert_0^2,$$
and thus $\Vert\cT\Vert_q\leq \Vert\cT\Vert_0$.
Since we can estimate
in the same way, by replacing $P_q$ by
$P_q^{-1}$
and $P_q'$ by
$P_q^{'-1}$, also
$\Vert\cT\Vert_0\leq\Vert\cT\Vert_q$, we get the assertion. \qed
\enddemo
\remark{1.5. Notation}
For a contraction $T:\HH\to\HH'$, we denote the extension
of $\cF(T)$ from $\cFf\to\cFfs$ to $\cFH\to\cFHs$ by $\cFT$.
\endremark
\remark{1.6. Remarks}
1) One might call $\cFT$ the second quantization of $T$, but we
will reserve this name for the restriction of $\cFT$ to some
operator algebra lying in $\cFH$ -- see the next section, where
we will also prove some positivity properties of this restricted
version. \newline
2) The operator $\cFT$ and its differential version (in particular
the number operator) were also considered in \cite{Wer} and
\cite{Sta,Mol}, respectively.\newline
3) It is clear that $\cFq(\cdot)$
behaves nicely with respect to composition and adjungation,
i.e.
$$\cFq(\id)=\id,\qquad
\cFq(ST)=\cFq(S)\cFT,\qquad\cFq(T^*)=\cFT^*,$$
but not with respect to the additive structure, i.e.
$$\cFq(T+S)\not=\cFT+\cFq(S)\qquad\text{in general.}$$
\endremark
In the context of the $q$-relations one usually encounters some 
kind of $q$-combi\-nato\-rics. Let us just remind of the basic facts.
\remark{1.7. Notations}
We put for $n\in\Nn$
$$\lb n\rb_q:=\frac{1-q^n}{1-q}=1+q+\dots+q^{n-1}
\qquad (\lb 0\rb_q:=0).$$
Then we have the $q$-factorial
$$\lb n\rb_q!:=\lb 1\rb_q\dots\lb n\rb_q,\qquad
\lb 0\rb_q!:=1$$
and a $q$-binomial coefficient
$$\binom nk_q:=\frac{\lb n\rb_q!}{\lb k\rb_q! \lb n-k\rb_q!}
=\prod_{i=1}^{n-k}\frac{1-q^{k+i}}{1-q^i}.$$
Another quite frequently used symbol is the $q$-analogue of the
Pochhammer symbol
$$(a;q)_n:=\prod_{j=0}^{n-1}(1-aq^j)\qquad\text{in particular}
\qquad (a;q)_\infty:=\prod_{j=0}^\infty(1-aq^j).$$
\endremark
The importance of these concepts in connection with the
$q$-relations can be seen from the following $q$-binomial
theorem, which is by now quite standard.
\proclaim{1.8. Proposition}
Let $x$ and $y$ be indeterminates which $q$-commute in the
sense $xy=qyx$. Then one has for $n\in\NN$
$$(x+y)^n=\sum_{k=0}^n\binom nk_q y^kx^{n-k}.$$
\endproclaim
\demo{Proof}
This is just induction and the easily checked equality
$$\binom nk_q+q^k\binom n{k+1}_q=\binom{n+1}{k+1}_q. \qed $$
\enddemo
In the same way as the usual Hermite polynomials are connected
to the bosonic relations, the $q$-relations are linked to
$q$-analogues of the Hermite polynomials.
\definition{1.9. Definition}
The polynomials $\Hnq$ ($n\in\Nn$), determined by
$$H_0^{(q)}(x)=1,\qquad H_1^{(q)}(x)=x$$
and
$$x\Hnq(x)=H_{n+1}^{(q)}(x)+\lb n\rb_q H_{n-1}^{(q)}(x)\qquad
(n\geq 1)$$
are called {\it $q$-Hermite polynomials}.
\enddefinition
We recall two basic facts about these polynomials which will
be fundamental for our investigations on 
the classical aspects of $q$-Gaussian processes.
\proclaim{1.10. Theorem}
1) Let $\nu_q$ be the measure on the interval
$\lb -2/\sqrt{1-q},2/\sqrt{1-q}\rb$ 
given by
$$\nu_q(dx)=\frac 1\pi \sqrt{1-q}\sin\theta\prod_{n=1}^\infty
(1-q^n)\vert 1-q^n e^{2i\theta}\vert^2 dx,$$
where
$$x=\frac 2{\sqrt{1-q}}\cos\theta\qquad\text{with $\theta\in
\lb 0,\pi\rb.$}$$
Then the $q$-Hermite polynomials are orthogonal with respect
to $\nu_q$, i.e.
$$\int_{-2/{\sqrt{1-q}}}^{2/{\sqrt{1-q}}} H_n(x)H_m(x)
\nu_q(dx)=\delta_{nm}\lb n\rb_q!.$$
2) Let $r>0$ and $x,y\in\lb -2/\sqrt{1-q},2/\sqrt{1-q}\rb$.
Denote by $p_r^{(q)}(x,y)$ the kernel
$$p_r^{(q)}(x,y):=\sum_{n=0}^\infty 
\frac{r^n}{\lb n\rb_q!}\Hnq(x)\Hnq(y).$$
Then we have with
$$x=\frac 2{\sqrt{1-q}}\cos\ff,\qquad y=\frac 2{\sqrt{1-q}}\cos\psi$$
the formula
$$p_r^{(q)}(x,y)=\frac{(r^2;q)_\infty}
{\vert(r e^{i(\ff+\psi)};q)_\infty
(re^{i(\ff-\psi)};q)_\infty\vert^2}.$$
In particular, for $q=0$, we get
$$p_r^{(0)}(x,y)=\frac{1-r^2}{(1-r^2)^2-r(1+r^2)xy+r^2(x^2+y^2)}.$$
\endproclaim
As usually in $q$-mathematics these formulas are quite old, 
namely the orthogonalizing measure $\nu_q$ was calculated by
Szego \cite{Sze}, whereas the kernel $p_r^{(q)}(x,y)$ goes even
back to Rogers \cite{Rog}. For more recent treatments, see
\cite{Bre,ISV,GR,LM}.
\bigskip
\heading
{\bf 2. Second quantization -- the functor $\sqq$}
\endheading
An abstract way of dealing with classical Gaussian processes is
by using the Gaussian functor $\Gamma$. This is a functor
from real Hilbert spaces and contractions to commutative von
Neumann algebras with specified trace-state and unital trace
preserving completely positive maps 
\cite{Nel1,Nel2,Gro,Sim1,Sim2}. A fermionic analogue of this
functor is also known, see, e.g., \cite{Wil,CL}.
\par
In this section we will present a $q$-analogue of the Gaussian
functor. Namely, to each real Hilbert space, $\HR$, we will
associate a von Neumann algebra with specified trace-state,
$(\sqH,\ee)$, and to every contraction $T:\HR\to\HR'$ a unital
completely positive trace preserving map $\sqT:\sqH\to
\sqq(\HR')$. 
\definition{2.1. Definition}
Let $\HR$ be a real Hilbert space and 
$\HC$ its complexification $\HC=\HR\oplus i\HR$.
Put, for $f\in\HR$,
$$\oo(f):=a(f)+a^*(f)\in B(\cFq(\HC))$$
and denote by $\sqH\subset B(\cFq(\HC))$ the von Neumann
algebra generated by all $\oo(f)$,
$$\sqH:=\vN(a(f)+a^*(f)\mid f\in\HR).$$
\enddefinition
\remark{2.2. Notation}
We denote by 
$$\ee:\sqH\to\CC$$
the vacuum expectation state on $\sqH$ given by
$$\ee\lb X\rb:=\la \Omega,X\Omega\ra_q\qquad (X\in\sqH).$$
\endremark 
We remind of some basic facts about $\sqH$ in the
following proposition.
\proclaim{2.3. Proposition}
The vacuum $\Omega$ is a cyclic and separating trace-vector
for $\sqH$, hence the vacuum expectation $\ee$ is a faithful normal
trace on $\sqH$ and $\sqH$ is a finite von Neumann algebra in
standard form.
\endproclaim
\demo{Proof}
See Theorems 4.3 and 4.4 in \cite{BSp3}.
\enddemo
The first part of the proposition yields in particular that the
mapping
$$\align
\sqH&\to\cFq(\HC)\\
X&\mapsto X\Omega
\endalign$$
is injective, in this way we can identify each $X\in\sqH$ 
with some element of the $q$-Fock space $\cFq(\HC)$.
\remark{2.4. Notations}
1) Let us denote by
$$\LuH:=\sqH\Omega$$
the image of $\sqH$ under the mapping $X\mapsto X\Omega$.
\newline
2) We also put
$$\LzH:=\cFq(\HC).$$
\endremark
\definition{2.5. Definition}
Let $\symb:\LuH\to\sqH$ be the identification of
$\LuH$ with $\sqH$ given by the requirement
$$\symb(\xi)\Omega=\xi\qquad\text{for}\quad
\xi\in\LuH\subset\LzH=\cFq(\HC).$$
\enddefinition
\remark{2.6. Remarks}
1) Of course, not each element of the $q$-Fock space
comes from an $X\in\sqH$, but the main
relation for observing the cyclicity of $\Omega$, namely
$$f_1\otimes\dots\otimes f_n=\oo(f_1)\dots \oo(f_n)\Omega-\eta
\qquad\text{with}\quad \eta\in\bigoplus_{l=0}^{n-1}\HR^{\otimes l},$$
yields that we have at least 
$f_1\otimes\dots\otimes f_n\in\LuH$.
\newline
2) In a quantum field theoretic context 
\cite{Sim1,Sim2} the operator $\symb(\fen)$
would be called \lq Wick product' and denoted by
$$\symb(\fen)=:\oo(f_1)\dots\oo(f_n):.$$
3) In a quantum probabilistic context 
\cite{Par2,Mey} $\symb$ would correspond
to taking an iterated quantum stochastic integral: For
$\HC=L^2(\RR)$ and $\xi=f_1\odo f_n$ with
$\xi(t_1,\dots,t_n)=f_1(t_1)\dots f_n(t_n)$
one would denote
$$\symb(\xi)=\int \xi(t_1,\dots,t_n)
d\oo(t_1)\dots d\oo(t_n)$$
and call $\xi$ the \lq Maassen kernel' of $\symb(\xi)$.
\endremark
The explicit form of our Wick products is given in the
following proposition.
\proclaim{2.7. Proposition}
We have for $n\in \NN$ and $f_1,\dots,f_n\in\HR$ the normal
ordered representation
$$\multline
\symb(\fen)=\\=\sum\Sb k,l=0,\dots,n\\k+l=n\endSb
\thinspace\sum\Sb
I_1=\{i(1),\dots,i(k)\}\\
I_2=\{j(1),\dots,j(l)\}\\
\text{with}\\
I_1\cup I_2=\{1,\dots,n\}\\
I_1\cap I_2=\emptyset\endSb
a^*(f_{i(1)})\dots a^*(f_{i(k)})a(f_{j(1)})\dots a(f_{j(l)})
\cdot q^{i(I_1,I_2)},\endmultline$$
where
$$i(I_1,I_2):=\#\{(p,q)\mid 1\leq p\leq k, 1\leq q\leq l,
i(p)>j(q)\}.$$ 
\endproclaim
Denote by $X$ the right hand side of the above relation. It
is clear that $X\Omega=\fen$, the problem is to see that
$X$ can be expressed in terms of the $\oo$'s.
\demo{Proof}
Note that the formula is true for
$$\symb(f)=\oo(f)=a(f)+a^*(f)$$
and that the definition of $a^*(f)$ and of $a(f)$ gives
$$\symb(f\otimes\fen)=
\oo(f)\symb(\fen)-\sum_{i=1}^n q^{i-1}\la f,f_i\ra
\symb(f_1\odo\check f_i\otimes\dots\otimes f_n).$$
>From this the assertion follows by induction. \qed
\enddemo
Note that $\symb(\fen)$ is just given by multiplying
out $\oo(f_1)\dots\oo(f_n)$ and bring all appearing terms
with the help of the relation $aa^*=qa^*a$ into a normal
ordered form -- i.e. we throw away all normal ordered terms
in $\oo(f_1)\dots\oo(f_n)$ which have less than $n$ factors.
Thus, for the special case $f_1=\dots=f_n$, we are in the 
realm of the $q$-binomial theorem and we have the following
nice formula.
\proclaim{2.8. Corollary}
We have for $n\in\NN$ and $f\in\HR$
$$\symb(f^{\otimes n})=\sum_{k=0}^n \binom nk_q a^*(f)^ka(f)^{n-k}.$$
\endproclaim
Instead of writing $\symb(\ftn)$ in a normal ordered form we can
also express it in terms of $\oo(f)$ with the help of the
$q$-Hermite polynomials.
\proclaim{2.9. Proposition}
We have for $n\in\Nn$ and $f\in\HR$ with $\Vert f\Vert=1$ the
representation
$$\symb(f^{\otimes n})=\Hnq(\oo(f)).$$
\endproclaim
\demo{Proof}
This follows by the fact that the $\symb(f^{\otimes n})$ fulfill the
same recurrence relation as the $\Hnq(\oo(f))$, namely
$$\oo(f)\symb(f^{\otimes n})=\symb(f^{\otimes (n+1)})+
\lb n\rb_q\symb(f^{\otimes(n-1)})$$
and that we have the same initial conditions
$$\symb(f^{\otimes 0})=\id,\qquad \symb(f^{\otimes 1})=\oo(f).\qed$$
\enddemo
We know \cite{Voi,VDN}
that for $q=0$ the von Neumann algebra $\sq_0(\HR)$ is isomorphic
to the von Neumann algebra of the free group on $\dim\HR$ generators --
in particular, it is a non-injective II$_1$-factor for $\dim\HR
\geq 2$. We conjecture non-injectivity and factoriality 
in the case $\dim\HR\geq 2$ for arbitrary $q\in(-1,1)$,
but up to now
we can only show the following.
\proclaim{2.10. Theorem}
1) For $-1<q<1$ and $\dim\HR> {16}/{(1-\vert q\vert)^2}$ the
von Neumann algebra $\sqH$ is not injective.
\newline
2) If $-1<q<1$ and $\dim\HR=\infty$ then $\sqH$ is a II$_1$-factor.
\endproclaim
\demo{Proof}
1) This was shown in a more general context in Theorem 4.2 in
\cite{BSp3}.
\newline
2) Let $\{e_i\}_{i\in\NN}$
be an orthonormal basis of $\HR$.
Fix $n\in\Nn$ and $r(1),\dots,r(n)\in\NN$ and consider
the operator
$$X:=\symb(e_{r(1)}\otimes \dots\otimes e_{r(n)}).$$
(For $n=0$ this shall be understood as $X=\id$.)
We put
$$\phi_m(X):=\frac 1m\sum_{i=1}^m \oo(e_i)X\oo(e_i)
\qquad (m\in\NN)$$
and claim that
$\phi_m(X)$ converges for $m\to\infty$ weakly to
$\phi(X):=q^nX$.
Because of the $m$-independent estimate
$$\Vert \phi_m(X)\Vert_q\leq \Vert X\Vert_q\thinspace\Vert \oo(e_1)
\Vert_q^2$$
it suffices to show
$$\lmu\la\xi,\phi_m(X)\eta\ra_q=\la\xi,\phi(X)\eta\ra_q$$
for all $\xi,\eta\in\cFq(\HC)$ of the form
$$\xi=e_{a(1)}\odo e_{a(u)},\qquad \eta=e_{b(1)}\odo e_{b(v)}$$
with $u,v\in\Nn$, $a(1),\dots,a(u),b(1),\dots,b(v)\in\NN$ (for $u=0$ we
put $\xi=\Omega$). 
To see this, put
$$m_0:=\max\{a(1),\dots,a(u),b(1),\dots,b(v),r(1),\dots,r(n)\}.$$
Since $\vert\la \xi,\oo(e_i)X\oo(e_i)\eta\ra_q\vert\leq M$ for some
$M$ (independent of $i$), we have
$$\align
\lmu\la\xi,\phi_m(X)\eta\ra_q&=\lmu \frac 1m\sum_{i=m_0+1}^m
\la\xi,\oo(e_i)X\oo(e_i)\eta\ra_q\\&=
\lmu\frac 1m\sum_{i=m_0+1}^m\la \xi,a(e_i)\symb(e_{r(1)}\odo
e_{r(n)})a^*(e_i)\eta\ra_q.\endalign$$
By Prop. 2.7, $\symb(e_{r(1)}\odo e_{r(n)})$ is now a linear
combination of terms of the form $Y=Y_1Y_2$ with
$$Y_1=a^*(e_{r(i(1))})\dots a^*(e_{r(i(k))})\qquad\text{and}
\qquad
Y_2=a(e_{r(j(1))})\dots a(e_{r(j(l))})$$
with $k+l=n$. Each such term gives, for $i>m_0$, a contribution
$$\align
\la\xi,a(e_i)Ya^*(e_i)\eta\ra_q&=
\la\xi,a(e_i)Y_1Y_2a^*(e_i)\eta\ra_q\\
&=q^{k+l}\la \xi,Y_1a(e_i)a^*(e_i)Y_2\eta\ra_q\\
&=q^{n}\la\xi,Y_1(\id+qa^*(e_i)a(e_i))Y_2\eta\ra_q\\
&=q^n\la\xi,Y_1Y_2\eta\ra_q\\
&=q^n\la\xi,Y\eta\ra_q\endalign$$
and hence
$$\align
\lmu\la\xi,\phi_m(X)\eta\ra_q&=\lmu\frac 1m\sum_{i=m_0+1}^m
q^n\la\xi,\symb(e_{r(1)}\odo e_{r(n)})\eta\ra_q
=\la\xi,q^nX\eta\ra_q.\endalign$$ 
Thus we have shown
$$\wlim_{m\to\infty}\phi_m(X)=\phi(X).$$
Let now $\tr$ be a normalized normal trace on $\sqH$. Then
$$\align
\tr\lb \phi(X)\rb &=\lmu\tr\lb \phi_m(X)\rb\\
&=\lmu\frac 1m\sum_{i=1}^m\tr\lb \oo(e_i)X\oo(e_i)\rb \\
&=\lmu\frac 1m\sum_{i=1}^m\tr\lb X\oo(e_i)\oo(e_i)\rb \\
&=\tr\lb X\cdot\lmu\frac 1m\sum_{i=1}^m\oo(e_i)\oo(e_i)\rb\\
&=\tr\lb X\phi(\id)\rb \\
&=\tr\lb X\rb .\endalign$$
Since $\phi^k(X)=q^{kn}X$ converges, for $k\to\infty$, (even in norm)
to 
$$\ee\lb X\rb\cdot\id=\cases 0,& n\geq 1\\
X=\id,& n=0,\endcases$$
we obtain
$$\align
\tr\lb X\rb&=\lku\tr\lb \phi^k(X)\rb 
=\tr\lb \lku\phi^k(X)\rb
=\ee\lb X \rb\tr\lb \id\rb
=\ee\lb X\rb.\endalign$$
Thus $\tr$ coincides on all operators of the form
$$X=\symb(e_{r(1)}\odo e_{r(n)})\qquad(n\in\Nn,\thinspace
r(1),\dots,r(n)\in\NN)$$
with our canonical trace $\ee$. Since the set of finite linear
combinations of such operators $X$ is weakly dense in $\sqH$,
we get the uniqueness of a normalized normal trace on
$\sqH$, which implies that $\sqH$ is a factor. \qed
\enddemo
The second part of our $q$-Gaussian functor $\sqq$ assigns to 
each contraction $T:\HR\to\HR'$ a map
$\sqT:\sqq(\HR)\to\sqq(\HR')$.
The idea is to extend $\sqT\oo(f)=\oo(Tf)$ in a canonical
way to all of $\sqH$. In general, the $q$-relations prohibit
the extension as a homomorphism, i.e.
$$\sqT\oo(f_1)\dots\oo(f_n)\not=\oo(Tf_1)\dots\oo(Tf_n)\qquad
\text{in general.}$$
But what can be done is to demand the above relation for the 
normal ordered form, i.e.
$$\sqT\symb(\fen)=\symb(Tf_1\odo Tf_n)
=\symb(\cFT \fen),$$
or
$$(\sqT X)\Omega=\cFT(X\Omega).$$
Thus our second quantization $\sqT$ is the restriction of
$\cFT$ from $\cFH= \LzH$ to $\sqH\cong\LuH$ and the question on
the existence of $\sqT$ amounts to the problem whether
$\cFT(\LuH)\subset L_q^{\infty}(\HR')$.
We know that $\cFT$ can be defined for $T$ a contraction and we
will see in the next theorem that no extra condition is needed
to ensure its nice behaviour with respect to $L_q^{\infty}$.
The case $q=0$ is due to Voiculescu \cite{Voi,VDN}.
\proclaim{2.11. Theorem}
1) Let $T:\HR\to\HR'$ be a contraction between real Hilbert
spaces. There exists a unique map
$\sqT:\sqq(\HR)\to\sqq(\HR')$
such that
$$(\sqT X)\Omega=\cFT (X\Omega).$$
The map $\sqT$ is linear, bounded, completely positive, unital
and preserves the canonical trace $\ee$.\newline
2) If $T$ is isometric, then $\sqT$ is a faithful homomorphism, and
if $T$ is the orthogonal projection onto a subspace, then
$\sqT$ is a conditional expectation.
\endproclaim
\demo{Proof}
Uniqueness of $\sqT$ follows from the fact that $\Omega$ is
separating for $\sqq(\HR')$. To prove the existence and the
properties of $\sqT$ we notice that any contraction $T$ can
be factored 
\cite{Hal} as $T=POI$ where
\newline
-- $I:\HH\to\cK=\HR\oplus\HR'$ is an isometric embedding
\newline
-- $O:\cK\to\cK$ is orthogonal \newline
-- $P:\cK=\HR\oplus\HR'\to\HR'$ is an orthogonal
projection onto a subspace.
\newline
Thus if we prove our assertions for each of these three cases
then we will also get the general statement for
$\sqT=\sqq(P)\sqq(O)\sqq(I)$.
\newline
a) Let $I:\HH\to\cK=\HR\oplus\HR'$ be an 
isometric embedding and $Q:\cK\to\cK$ the orthogonal projection
onto $\HR$. Then $\cFq(Q)$ is a projection in $\cFq(\cKC)$ and
$\cFq(\HC)$ can be identified with $\cFq(Q)\cFq(\cKC)$. Let
us denote by $\oo_{\cK}(f)$ the sum of creation and annihilation
operator on $\cFq(\cKC)$. If we put
$$\sqq^\cK(\HR):=\vN(\oo_{\cK}(f)\mid f\in\HR)\subset
B(\cFq(\cKC)),$$
then
$$\sqq^\cK(\HR)\cFq(\HC)\subset\cFq(\HC)$$
and we have the canonical identification
$$\sqH\cong \sqq^\cK(\HR) \cFq(Q),$$
which gives a homomorphism (and thus a completely positive)
$$\sqq(I):\sqH\to\sqq(\cK).$$
Faithfulness is clear since $\cFq(Q)\Omega=\Omega$ and
$\Omega$ separating. This yields also that the trace is
preserved.
\newline
b) Let $P:\cK=\HR\oplus\HR'\to\HR'$ be an orthogonal
projection, i.e. $PP^*=\id_{\HR'}$, where $P^*:\HR'\to\cK$ is
the canonical inclusion. Then
$$\sqq(P)X:=\cFq(P)X\cFq(P^*)\qquad (X\in\sqq(\cK))$$
gives the right operator, because we have for $k,l\in\Nn$ and
$f_1,\dots,f_k,g_1,\dots,g_l\in\cK$
$$\align
\cFq(P)a^*(f_1)\dots a^*(f_k)&a(g_1)\dots a(g_l)\cFq(P^*)=\\
&=a^*(Pf_1)\dots a^*(Pf_k)\cFq(P)\cFq(P^*)
a(Pg_1)\dots a(Pg_l)\\
&=a^*(Pf_1)\dots a^*(Pf_k)a(Pg_1)\dots a(Pg_l).
\endalign$$
By its concrete form, $\sqq(P)$ is a conditional expectation
and
$$\ee\lb \cFq(P)X\cFq(P^*)\rb=
\la\cFq(P^*)\Omega,X\cFq(P^*)\Omega\ra_q=
\la\Omega,X\Omega\ra_q=\ee\lb X\rb$$
shows that it preserves the trace.
\newline
c) Let $O:\cK\to\cK$ be orthogonal, i.e. $OO^*=O^*O
=\id_\cK$.
Then, as in b),
$$\sqq(O)X=\cFq(O)X\cFq(O^*),$$
which is, by
$$\cFq(O^*)\cFq(O)=\cFq(\id_\cK)=\id_{\cFq(\cK_{\CC})}$$
also a faithful homomorphism. \qed
\enddemo
Instead of working on the level of von Neumann algebras we could
also consider the $C^*$-analogues of the above constructions.
This would be quite similar. We just indicate the main points.
\definition{2.12. Definition}
Let $\HR$ be a real Hilbert space and 
$\HC$ its complexification $\HC=\HR\oplus i\HR$.
Put, for $f\in\HR$,
$$\oo(f):=a(f)+a^*(f)\in B(\cFq(\HC))$$
and denote by $\PhqH\subset B(\cFq(\HC))$ the $C^*$-algebra
generated by all $\oo(f)$,
$$\PhqH:=C^*(a(f)+a^*(f)\mid f\in\HR).$$
\enddefinition
Clearly, the vacuum is also a separating trace-vector for 
$\PhqH$ and, by Remark 2.6., it is also cyclic and
$\symb(\fen)\in\PhqH$ for all $n\in\Nn$ and all $f_1,\dots,
f_n\in\HR$.
\par
The most important fact for our latter considerations is that
$\sqq(T)$ can also be restricted to the $C^*$-level.
\proclaim{2.13. Theorem}
1) Let $T:\HR\to\HR'$ be a contraction between real Hilbert
spaces. There exists a unique map
$\PhqT:\Phq(\HR)\to\Phq(\HR')$
such that
$$(\PhqT X)\Omega=\cFT (X\Omega).$$
The map $\PhqT$ is linear, bounded, completely positive, unital
and preserves the canonical trace $\ee$.\newline
2) If $T$ is isometric, then $\PhqT$ is a faithful homomorphism, and
if $T$ is the orthogonal projection onto a subspace, then
$\PhqT$ is a conditional expectation.
\newline
3) We have $\PhqT=\sqT/\PhqH$.
\endproclaim
\demo{Proof}
This is analogous to the proof of Theorem 2.11. \qed
\enddemo
We can now also prove the analogue of the second part of Theorem
2.10. The analogue 
of factoriality for $C^*$-algebras is simplicity.
\proclaim{2.14. Theorem}
If $-1<q<1$ and $\dim\HR=\infty$ then $\PhqH$ is simple.
\endproclaim
\demo{Proof}
Again, this is similar to the proof of the von Neumann algebra result.
We just indicate the main steps. \newline
We use the notations from the proof of Theorem 2.10. First, by
norm estimates, one can show that
the convergence
$\lmu\phi_m(X)=\phi(X)$ for $X$ of the form
$X:=\symb(e_{r(1)}\otimes \dots\otimes e_{r(n)})$
is even a convergence in norm. Since $\phi(X)$ is nothing but
$\phi(X)=\sqq(q)X$, where $q$ is regarded as multiplication 
operator on $\HR$, we have, by 2.13, the bound
$$\Vert\phi(X)\Vert_q\leq\Vert X\Vert_q.$$
This together with the $m$-independent bound
$$\Vert\phi_m(X)\Vert_q\leq\Vert X\Vert_q\thinspace\Vert\oo(e_1)
\Vert_q^2$$
implies that
$$\lmu\phi_m(X)=\sqq(q)X\qquad\text{uniformly for all $X\in\PhqH$.}$$
Now assume we have a non-trivial ideal $I$ in $\PhqH$ and consider
a positive non-vanishing $X\in I$. Then $\phi_m(X)\in I$ for all
$m\in \NN$ and thus $\sqq(q)X\in I$. Iterating shows $\sqq(q^n)X\in I$
for all $n\in \NN$ and because of the uniform convergence
$\lim_{n\to\infty}\sqq(q^n)X=\ee\lb X\rb\id$
we obtain $\ee\lb X\rb\id\in I$. The faithfulness of $\ee$
implies then
$I=\PhqH$. \qed
\enddemo
\remark{2.15. Remark}
One might be tempted to conjecture that, for fixed $\HR$, the 
von Neumann algebras $\sqH$ or the $C^*$-algebras $\PhqH$ are for
all $q\in(-1,1)$ isomorphic. At the moment, no results in this
direction are known. One should note that there exist partial
answers \cite{JSW1,JSW2,JW,DN} to the analogous question for the
$C^*$-algebra generated by $a(f),a^*(f)$ (not the sum) showing
that at least for small values of $q$ 
and $n:=\dim\HC<\infty$ these algebras are isomorphic to the
$(q=0)$-algebra, which is
an extension of the Cuntz algebra $O_n$ by compact operators
\cite{Cun,Eva}.
However, the methods used there do not extend to the case of
$\sqH$ or $\PhqH$.
\endremark
\bigskip
\heading
{\bf 3. Non-commutative aspects of $q$-Gaussian processes}
\endheading
Before we define the notion of a $q$-Gaussian process, we want to
present our general frame on non-commutative processes.
By $T$ we will denote the range of our time parameter $t$, typically
$T$ will be some interval in $\RR$.
\definition{3.1. Definitions}
1) Let $\cA$ be a finite von Neumann algebra and $\ff:\cA\to\CC$
a faithful normal trace on $\cA$. Then we call the pair
$(\cA,\ff)$ a {\it (tracial) probability space}. \newline
2) A {\it random variable} on $(\cA,\ff)$ is a 
self-adjoint operator $X\in\cA$.
\newline
3) A {\it stochastic process} on $(\cA,\ff)$ is a family $(X_t)_{t\in T}$
of random variables $X_t\in\cA$ ($t\in T$).
\newline
4) The {\it distribution} of a random variable $X$ on $\Aff$ is the
probability measure $\nu$ on the spectrum of $X$ determined by
$$\ff(X^n)=\int x^n d\nu(x)\qquad\text{for all $n\in\Nn$.}$$
\enddefinition
We should point out that there are also a lot of 
quantum probabilistic investigations
in the context of more general, non-tracial situations, see e.g.
\cite{AFL,K\"um1}. Of course,
life becomes much harder there.
\par
We will only consider centered Gaussian processes, thus a $q$-Gaussian
process will be totally determined by its covariance.
Since we would like 
to have realizations of our processes on {\it separable}
Hilbert spaces, our admissible covariances are not just positive
definite functions, but they should admit a separable representation.
\definition{3.2. Definition}
A function $c:T\times T\to\RR$ is called {\it covariance function}, if
there exists a separable real Hilbert space $\HR$ and vectors
$f_t\in\HR$ for all $t\in T$ such that
$$c(s,t)=\la f_s,f_t\ra\qquad (s,t\in\HR).$$
\enddefinition
\definition{3.3. Definition}
Let $c:T\times T\to \RR$ be a covariance function corresponding
to a real Hilbert space
$\HR$ and vectors $f_t\in\HR$ ($t\in T$). Then we put for all $t\in T$
$$X_t:=\oo(f_t)\in\sqH$$
and call the process $\Xt$ on $(\sqH,\ee)$ the
{\it $q$-Gaussian process with covariance $c$}.
\enddefinition
\remark{3.4. Remarks}
1) Of course, the $q$-Gaussian process depends, up to isomorphism,
only on $c$ and not on the special choice of $\HH$ and $(f_t)_{t\in T}$.
\newline
2) We can characterize our $q$-Gaussian process also by the $q$-relations
in the form
$$X_t=a_t+a_t^*\qquad\text{and}\qquad \ee\lb
\thinspace\cdot\thinspace\rb=\la\Omega,\cdot
\thinspace\Omega\ra,$$
where for all $s,t\in T$
$$a_sa_t^*-qa_t^*a_s=c(s,t)\cdot\id\qquad\text{and}\qquad a_t\Omega=0.$$ 
In this form our $q$-Gaussian processes were introduced by Frisch
and Bourret \cite{FB}.
\endremark
We can now define $q$-analogues of all classical Gaussian processes,
just by choosing the appropriate covariance. In the following we
consider three prominent examples.
\definition{3.5. Definitions}
1) The $q$-Gaussian process $\XBM$ with covariance
$$c(s,t)=\min(s,t)\qquad (0\leq s,t<\infty)$$
is called {\it $q$-Brownian motion}. \newline
2) The $q$-Gaussian process $\XBB$ with covariance
$$c(s,t)=s(1-t)\qquad (0\leq s\leq t\leq 1)$$
is called {\it $q$-Brownian bridge}. \newline
3) The $q$-Gaussian process $\XOU$ with covariance
$$c(s,t)=e^{-\vert t-s\vert}\qquad (s,t\in\RR)$$
is called {\it $q$-Ornstein-Uhlenbeck process}.
\enddefinition
\remark{3.6. Remarks}
1) That the three examples for $c$ are indeed covariance functions
is clear by the existence of the respective classical processes,
for a direct proof see, e.g., \cite{Sim2}. \newline
2) The Ornstein-Uhlenbeck process is often also called 
oscillator process, see \cite{Sim2}.
\endremark
Let $(\cA,\ff)$ be a tracial probability space and let $\cB$ be
a von Neumann subalgebra of $\cA$. Then we have (see, e.g., 
\cite{Tak}) 
a unique conditional
expectation (\lq partial trace')
from $\cA$ onto $\cB$ which preserves the trace $\ff$ --  
which we will denote in a probabalistic
language by $\ff\lb\thinspace\cdot\thinspace\vert\cB\rb$.
Thus in the frame of tracial probability spaces we always have the
following canonical generalization of the classical Markov property
(which says that the future depends on the past only through the 
present).
\definition{3.7. Definition}
Let $(\cA,\ff)$ be a probability space and $\Xt$ a stochastic
process on $\Aff$. Denote by
$$\align
\cAbt&:=\vN(X_u\mid u\leq t)\subset\cA\\ 
\cAgt&:=\vN(X_t)\qquad\quad\subset\cA.\endalign$$
We say that $\Xt$ is a {\it Markov process} if we have for 
all $s,t\in T$ with $s\leq t$ the property
$$\ff\lb X\vert\cAbs\rb\subset\cAgs\qquad\text{for all
$X\in\cAgt$.}$$
\enddefinition
Now, the conditional expectations 
$\ee\lb\thinspace\cdot\thinspace\vert\cAbs\rb$ 
in the case of $q$-Gaussian
processes are quite easy to 
handle because they are nothing but the
second quantization of projections in the underlying Hilbert
space. Namely, consider a $q$-Gaussian process $\Xt$ corresponding
to the real Hilbert space $\HR$ and vectors $f_t$ ($t\in T$). Let
us denote by 
$$\align
\cHbt&:=\Span(f_u\mid u\leq t)\subset\HR\\
\cHgt&:=\RR f_t\quad\qquad\qquad\subset\HR\endalign$$
the Hilbert space analogues of $\cAbt$ and $\cAgt$, respectively.
Then we have
$$\cAbt\cong\sqq(\cHbt)\qquad\text{and}\qquad
\cAgt\cong\sqq(\cHgt),$$
and $\ee\lb\thinspace\cdot\thinspace
\vert\cAbt\rb=\sqq(\cPbt)$ is the second
quantization of the orthogonal projection
$$\cPbt:\HR\to\cHbt.$$
Thus we can translate the Markov property for $q$-Gaussian
processes into the following Hilbert space level statement.
\proclaim{3.8. Proposition}
Let $\Xt$ be a $q$-Gaussian process as above. It has the Markov
property if and only if
$$\cPbs\cHgt\subset\cHgs \qquad\text{for all $s,t\in T$ with
$s\leq t$.}$$
\endproclaim
Thus Markovianity is a property of the underlying Hilbert space and
does not depend on $q$ and we get as in the classical case the
following characterization in terms of the covariance.
\proclaim{3.9. Proposition}
A $q$-Gaussian process with covariance $c$ is Markovian if and
only if we have for all triples $s,u,t\in T$ with $s\leq u\leq t$
that
$$c(t,s)c(u,u)=c(t,u)c(u,s).$$
\endproclaim
\demo{Proof}
See the proof of Theorem 3.9 in \cite{Sim2}. \qed
\enddemo
\proclaim{3.10. Corollary}
The $q$-Brownian motion $\XBM$, the $q$-Brownian bridge $\XBB$, and
the $q$-Ornstein-Uhlenbeck process $\XOU$ are all Markovian.
\endproclaim
Analogously, we have all statements of the classical Gaussian
processes which depend only on Hilbert space properties. Let us
just state the characterization of the Ornstein-Uhlenbeck process
as the only stationary 
Gaussian Markov process 
with continuous covariance and the characterization
of martingales among the Gaussian processes.
\proclaim{3.11. Proposition}
Let $\Xt$ be a $q$-Gaussian process which is stationary, Markovian
and whose covariance $c(s,t)=c'(t-s)$ is continuous. Then
$X_t=\alpha X^{qOU}_{\beta t}$
for suitable $\alpha,\beta>0$.
\endproclaim
\demo{Proof}
See the proof of the analogous statement for classical Gaussian
processes, Corollary 4.10 in \cite{Sim2}. \qed
\enddemo
\definition{3.12. Definition}
Let $\Xt$ be a stochastic process on a probability
space $\Aff$ and let the notations be as in Definition 3.7.
Then we say that $\Xt$ is a {\it martingale} if 
$$\ff\lb X_t\vert\cAbs\rb= X_s\qquad\text{for all $s\leq t$.}$$
\enddefinition
\proclaim{3.13. Proposition}
A $q$-Gaussian process is a martingale if and only if
$\cPbs f_t=f_s$ for all $s\leq t$ -- which is the case if and
only if 
$c(s,t)=c(s,s)$ for all $s\leq t$.
\endproclaim
\demo{Proof}
We have
$$\oo(f_s)=X_s=\ee\lb X_t\vert\cAbs\rb=\sqq(\cPbs)\oo(f_t)=
\oo(\cPbs f_t),$$
implying $\cPbs f_t=f_s$. \qed
\enddemo  
\bigskip
\heading
{\bf 4. Classical aspects of $q$-Gaussian processes}
\endheading
In this section we want to address the question whether our
non-commutative stochastic processes can also be interpreted
classically.
\definition{4.1. Definition}
Let $\Xt$ be a stochastic process on some non-commutative
probability space $\Aff$. We call a classical
real-valued process $\Xtk$ on some classical probability space
$(\Omega,{\frak A},P)$ a {\it classical version} of
$\Xt$ if all time-ordered moments of $\Xt$ and $\Xtk$ coincide,
i.e. if we have for all $n\in\NN$, all $t_1\dots,t_n\in T$ with
$t_1\leq \dots\leq t_n$, and all bounded Borel functions $h_1,
\dots,h_n$ on $\RR$ the equality
$$\ff\lb h_1(X_{t_1})\dots h_n(X_{t_n})\rb=
\int_\Omega h_1(\tilde X_{t_1}(\omega))\dots h_n(
\tilde X_{t_n}(\omega))dP
(\omega).$$
\enddefinition
\remark{4.2. Remark}
Most calculations in a non-commutative context involve only
time-ordered moments, see, e.g., the calculation of the
Green function of the $q$-Ornstein-Uhlenbeck process in
\cite{NSp}. Thus, results of such calculations can also be interpreted
as results for the classical version -- if such a version exists.
\endremark
It is clear that there is at most one classical version for a 
given non-commutative process $\Xt$. The problem consists in
showing the existence.
If we denote by $\id_B$ the
characteristic function of a measurable subset $B$ of $\RR$, then we 
can construct the classical version $\Xtk$ of $\Xt$ via
Kolmogorov's existence theorem from the collection of all
$\mu_{t_1,\dots,t_n}$ ($n\in\NN$, $t_1\leq\dots\leq t_n$) -- which
are for $B_1,\dots,B_n\subset\RR$ defined by
$$\align
\mu_{t_1,\dots,t_n}(B_1\times\dots\times B_n)&=
P(\tilde X_{t_1}\in B_1,\dots,\tilde X_{t_n}\in B_n)\\
&=\ff\lb \id_{B_1}(X_{t_1})\dots \id_{B_n}(X_{t_n})\rb\endalign$$
-- if and only if all $\mu_{t_1,\dots,t_n}$ are probability measures.
Whereas this is of course the case for $\mu_{t_1}$ and, 
in our tracial frame because of
$$\mu_{t_1,t_2}(B_1\times B_2)=\ff\lb \id_{B_1}(X_{t_1})
\id_{B_2}(X_{t_2})\rb=
\ff\lb \id_{B_1}(X_{t_1})
\id_{B_2}(X_{t_2})\id_{B_1}(X_{t_1})\rb,$$
also for $\mu_{t_1,t_2}$,
there is no apriori reason why it should
be true for bigger $n$. And in general it is not. It
is essentially the content of Bell's inequality that there are
examples of non-commutative processes which possess no classical
version -- for a discussion of these subjects see, e.g., \cite{KM}.
\par
But for special classes of non-commutative processes classical versions
might exist. One prominent example of such a class are the Markov
processes.
\par
\definition{4.3. Definition}
Let $\Xt$ be a Markov process on a probability space $\Aff$.
Let, for $t\in T$, 
$\spec(X_t)$ and  $\nu_t$ be the spectrum und the distribution,
respectively, of the self-adjoint operator $X_t$. Denote by
$$L^\infty(X_t):=\vN(X_t)=L^\infty(\spec(X_t),\nu_t).$$
The operators
$$\Kst:L^\infty(X_t)\to L^\infty(X_s)\qquad (s\leq t),$$
determined by
$$\ff\lb h(X_t)\vert\cAbs\rb=
\ff\lb h(X_t)\vert \cAgs\rb=(\Kst h)(X_s)$$
are called {\it transition operators} of the process $\Xt$,
and, looked upon from the other side,
the process $\Xt$ is called a {\it dilation} of the
transistion operators $\cK=(\Kst)_{s\leq t}$.
\enddefinition
The following theorem is by now some kind of folklore in
quantum probability, see, e.g. \cite{AFL,K\"um2,BP,
Bia1}. We just indicate the
proof for sake of completeness.
\proclaim{4.4. Theorem} If $\Xt$ is a Markov process on some
probability space $\Aff$, then there exists a classical version
$\Xtk$ of $\Xt$, which is a classical Markov process.
\endproclaim
\demo{Proof}
One can express the time-ordered moments of a Markov process
in terms of the transition operators via
$$\align
\ff\lb h_1(X_{t_1})\dots h_n(X_{t_n})\rb&=
\ff\lb h_1(X_{t_1})\dots h_n(X_{t_n})\vert \cA_{t_{n-1}\rb}\rb\\
&=\ff\lb h_1(X_{t_1})\dots h_{n-1}(X_{t_{n-1}})\ff\lb h_n(X_{t_n})
\vert\cA_{t_{n-1}\rb}\rb\rb\\
&=\ff\lb h_1(X_{t_1})\dots h_{n-1}(X_{t_{n-1}})(\cK_{t_{n-1},t_n}
h_n)(X_{t_{n-1}})\rb\\
&=\ff\lb h_1(X_{t_1})\dots h_{n-2}(X_{t_{n-2}})(h_{n-1}\cdot
\cK_{t_{n-1},t_n}h_n)
(X_{t_{n-1}})\rb\\
&=\dots\\
&=\ff\lb (h_1\cdot \cK_{t_1,t_2}(h_2\cdot\cK_{t_2,t_3}(h_3\cdot\dots)))
(X_{t_1})\rb,\endalign$$
from which it follows -- because $\Kst$ preserves positivity --
that the corresponding $\mu_{t_1,\dots,t_n}$ are probability 
measures. 
That the classical version is also a classical Markov process
follows by the same formula.
\qed
\enddemo
\proclaim{4.5. Corollary}
There exist classical versions of all $q$-Gaussian Markov processes.
In particular, we have classical versions of 
the $q$-Brownian motion,
of the $q$-Brownian bridge, and of the
$q$-Ornstein-Uhlenbeck process.
\endproclaim
Our aim now is to describe these classical versions more
explicitly by calculating their transition probabilities in
terms of the orthogonalizing measure $\nu_q$ and the kernel
$p_r^{(q)}(x,y)$ of Theorem 1.10.
\proclaim{4.6. Theorem}
Let $\Xt$ be a $q$-Gaussian Markov process with covariance $c$
and put
$$\lt:=\sqrt{c(t,t)}\qquad\text{and}\qquad\lst:=\frac{c(t,s)}
{\sqrt{c(s,s)c(t,t)}}.$$
1) We have
$$L^\infty(X_t)=L^\infty(\lb -{2\lt}/{\sqrt{1-q}},
{2\lt}/{\sqrt{1-q}}\rb,\nu_q(dx/\lt)).$$
2) If $\lst=\pm1$, then the transition operator $\Kstq$ is given by
$$(\Kstq h)(x)=h(\pm x\lt/\ls).$$
If $\vert \lst\vert<1$, then
the transition operator $\Kstq$ is given by
$$(\Kstq h)(x)=\int h(y)\kstq(x,dy),$$
where the transition probabilities $\kstq$ are Feller kernels 
which have the explicit form
$$\kstq(x,dy)=p_\lst^{(q)}(x/\ls,y/\lt)\nu_q(
{dy}/\lt).$$
In particular, for $q=0$ and $\vert\lst\vert<1$, 
we have the following transition
probabilities for the free Gaussian Markov processes
$$\multline
\kstn(x,dy)=\\=\frac 1{2\pi\ltq}\frac
{(1-\lstq)\sqrt{4\ltq-y^2}dy}
{(1-\lstq)^2-\lst(1+\lstq) (x/\ls) (y/\lt)+
\lstq(({x^2}/{\lsq})+({y^2}/{\ltq}))}.\endmultline$$
\endproclaim
Recall that a kernel $k(x,dy)$ is called Feller, if the map
$x\mapsto k(x,dy)$ is weakly continuous and 
$k(x,\cdot)\to 0$ weakly as $x\to\pm\infty$ -- or equivalently that
the corresponding 
operator $\cK$ sends $C_0(\RR)$ to $C_0(\RR)$, see, e.g.,
\cite{DM}.
\demo{Proof}
1) This was shown in \cite{BSp2}; noticing the connection 
between $q$-relations and $q$-Hermite polynomials the assertion
reduces essentially to part 1) of Theorem 1.10. 
\newline
2) By Prop. 2.9, we know
$$\symb(\ftn)=\Vert f\Vert^n \Hnq({\oo(f)}/{\Vert f\Vert}).$$
Let our $q$-Gaussian process $\Xt$ now be of the form
$X_t=\oo(f_t)$. Markovianity implies
$$\cPbs f_t=\mu f_s\qquad\text{where}\qquad
\mu=\frac{\la f_t,f_s\ra}{\la f_s,f_s\ra}=\frac {c(t,s)}{c(s,s)}.$$
Because of
$$\ee\lb\symb(\ftn_t)\vert\cAbs\rb=\symb((\cPbs f_t)^{\otimes n})
=\mu^n\symb(\ftn_s)$$
we obtain with
$$\lt:=\Vert f_t\Vert=\sqrt{c(t,t)}\qquad\text{and}\qquad
\lst:=\mu\frac \ls\lt=\frac{c(t,s)}{\sqrt{c(s,s)c(t,t)}}$$
the formula
$$\align
\ee\lb\Hnq({X_t}/\lt)\vert\cAbs\rb&=
\frac 1{\ltn}\ee\lb\symb(\ftn_t)\vert\cAbs\rb\\
&=\frac{\mu^n}{\ltn}\symb(\ftn_s)\\
&=(\mu\frac \ls\lt)^n\Hnq({X_s}/\ls)\\
&=\lstn\Hnq({X_s}/\ls),
\endalign$$
implying
$$\Kstq(\Hnq(\cdot/\lt))=\lstn\Hnq(\cdot/\ls).$$
Let us now consider the canonical 
extension of our transition operators from the $L^\infty$-spaces to the 
$L^2$-spaces, i.e.
$$\Kstq:L^2(X_t)\to L^2(X_s).$$
If we use the fact that the rescaled $q$-Hermite polynomials
$(\Hnq(\cdot/\lt)/\sqrt{\lb n\rb!})_{n\in\Nn}$ constitute an orthonormal
basis of $L^2(X_t)$, we get directly the assertion
in the case $\lst=\pm 1$. (For $\lst=-1$ one also has to note that
$H^{(q)}_{2k}$ and $H^{(q)}_{2k+1}$ are even and odd polynomials,
respectively.) 
\newline
In the case $\vert \lst\vert<1$, our formula implies that $\Kstq$ is
a Hilbert-Schmidt operator, thus it has a concrete representation
by a kernel $\kstq$, which is given by
$$\align
\kstq(x,dy)&=\sum_{n=0}^\infty\frac{\lstn}{\lb n\rb_q!}
\Hnq(x/\ls)\Hnq(y/\lt)\nu_q({dy}/\lt)\\
&=p_\lst(x/\ls,y/\lt)\nu_q({dy}/\lt).\endalign$$
That our kernels are Feller follows from the fact that, by
Theorem 2.13, our second quantization (i.e. our transition
operators) restrict to the $C^*$-level (i.e. to continuous
functions).
\newline
The formula for $\kstn$ follows from the concrete form of
$p_r^{(0)}$ of Theorem 1.10 and the fact that
$$\nu_0(dy)=\frac 1{2\pi}\sqrt{4-y^2}dy\qquad\text{for}
\qquad y\in\lb -2,2\rb. \qed $$
\enddemo
The main formula of our proof, namely the action of the
conditional expectation on the $q$-Hermite polynomials, says that
we have some quite canonical martingales associated to $q$-Gaussian
Markov processes -- provided the factor $\lst$ decomposes into a quotient
$\lst=\lambda(s)/\lambda(t)$. Since this can be assured by a corresponding
factorization property of the covariance function -- which is not very
restrictive for Gaussian Markov processes, see Theorem 4.9 of
\cite{Sim2} -- we get the following corollary.
\proclaim{4.7. Corollary}
Let $\Xt$ be a $q$-Gaussian process whose covariance factorizes 
for suitable functions $g$ and $f$ as
$$c(s,t)=g(s)f(t) \qquad \text{for $s\leq t$}.$$
Then, for all $n\in\Nn$,  the processes $(M_n(t))_{t\in T}$ with
$$M_n(t):=\bigl({g(t)}/{f(t)}\bigr)^{n/2} \Hnq(X_t/\lt)$$
are martingales.
\endproclaim
Note that the assumption on the factorization of the covariance is
in particular fulfilled for the $q$-Brownian motion, for the
$q$-Ornstein-Uhlenbeck process, and for the $q$-Brownian bridge.
\demo{Proof}
Our assumption on the covariance implies
$$\lst=\sqrt{\frac{g(s)/f(s)}{g(t)/f(t)}},$$
hence our formula for the action of the conditional expectation
on the $q$-Hermite polynomials (in the proof of Theorem 4.6) can
be written as
$$\bigl({g(t)}/{f(t)}\bigr)^{n/2}\ee\lb\Hnq({X_t}/\lt)\vert\cAbs\rb
=\bigl({g(s)}/{f(s)}\bigr)^{n/2}\Hnq(X_s/\ls),$$
which is exactly our assertion. \qed
\enddemo
\remark{4.8. Remark}
Consider the $q$-Brownian motion $\XBM$. Then the Corollary states
that 
$$M_n^{(q)}(t):=t^{n/2}\Hnq(X_t^{qBM}/\sqrt t)$$
is a martingale. In terms of quantum stochastic integrals
these martingales would have the form
$$M_n^{(q)}(t)=\intdotsint\Sb {0\leq t_1,\dots,t_n\leq t}\\ t_i\not= t_j 
\thinspace (i\not= j)\endSb  dX_{t_1}^{qBM}\dots
dX_{t_n}^{qBM}.$$
Since at the moment, for general $q$, no
rigorous theory of $q$-stochastic integration exists, this has to
be taken as a purely formal statement. 
For $q=0$, however, such a rigorous theory was developed in 
\cite{KSp}, and the above representation by stochastic integrals was
established by Biane \cite{Bia2}. In this case, he could put
this representation into the form of the stochastic differential
equation
$$M_n^{(0)}(t)
=\sum_{k=0}^{n-1}\int_0^t M_k^{(0)}(s)dX_s^{0BM}M^{(0)}_{n-k-1}(s),$$
which should be compared with the classical formula
$$M_n^{(1)}(t)=n\int_0^t M^{(1)}_{n-1}(s) dX_s^{1BM}.$$
\endremark 
\example{4.9. Example: free Gaussian processes}
We will now specialize the formula for $\kstn$ to the case
of the free Brownian motion, the free Ornstein-Uhlenbeck process
and the free Brownian bridge. The transition probabilities for
the two former cases were also derived by Biane \cite{Bia1}
in the context of processes with free increments. 
1) free Brownian motion: We have $c(s,t)=
\min(s,t)$, thus 
$$\lt=\sqrt t\qquad\text{and}\qquad \lst=\sqrt{s/t}.$$  
This yields
$$\ck_{s,t}(x,dy)=
\frac
{(t-s)}
{(t-s)^2-(t+s)xy+x^2t+y^2s}
\frac{\sqrt{4t-y^2}dy}{2\pi}$$
for 
$$x\in\lb -2\sqrt s,2\sqrt s\rb\qquad\text{and}\qquad
y\in\lb -2\sqrt t,2\sqrt t\rb.$$
\newline
2) free Ornstein-Uhlenbeck process: We have
$c(s,t)=e^{-\vert t-s\vert}$, thus 
$$\lt=1\qquad\text{and}\qquad 
\lst=e^{-\vert t-s\vert}.$$ 
Since this process is stationary, it
suffices to consider the transition probabilities for $s=0$:
$$\ck_{0,t}(x,dy)=\frac
{(e^{2t}-1)}
{4\sinh^2 t-2xy\cosh t+x^2+y^2}\frac{\sqrt{4-y^2}dy}{2\pi}
\qquad\text{for}\qquad x,y\in \lb-2,2\rb.$$
Let us also calculate the generator $N$ of this process --
which is characterized by $$\Kst=e^{-(t-s)N}.$$ 
It has
the property
$$N H_n^{(0)}=nH_n^{(0)}\qquad (n\in\Nn),$$
and
differentiating the above kernel shows that it should be given
formally by a kernel $-2/(y-x)^2$ with respect to $\nu_0$. Making
this more rigorous \cite{vWa} yields that $N$ has on functions
which are differentiable the form
$$(Nh)(x)=xf'(x)-2 \int\frac{f(y)-f(x)-f'(x)(y-x)}
{(y-x)^2}\nu_0(dy).$$
3) free Brownian bridge: We have 
$c(s,t)=s(1-t)$ for $s\leq t$, thus 
$$\lt=\sqrt{t(1-t)}\qquad\text{and}\qquad
\lst=\sqrt{\frac{s(1-t)}{t(1-s)}}.$$ 
This yields
$$\multline\ck_{s,t}(x,dy)=\\=
\frac {1-s}{1-t}\frac{(t-s)}
{(t-s)^2-(s+t-2st)xy+t(1-t)x^2+s(1-s)y^2}
\frac{\sqrt{4t(1-t)-y^2}dy}{2\pi},\endmultline$$
for
$$x\in\lb-2\sqrt{s(1-s)},2\sqrt{s(1-s)}\rb\qquad\text{and}\qquad
y\in\lb-2\sqrt{t(1-t)},2\sqrt{t(1-t)}\rb.$$
\endexample
\example{4.10. Example: fermionic Gaussian processes}
For illustration, we also want to consider the fermionic
($q=-1$) analogue of Gaussian processes. Although this case
has not been included in our frame everything works similar,
the only difference is that in the Fock space we get a kernel
of our scalar product consisting of anti-symmetric tensors. This
is responsible for the fact that the corresponding ($-1$)-Hermite
polynomials collapse just to 
$$H_0^{(-1)}(x)=1\qquad\text{and}\qquad
H_1^{(-1)}(x)=x.$$
The corresponding measure $\nu_{-1}$ is not absolutely continuous
with respect to the Lebesgue measure anymore, but collapses to 
$$\nu_{-1}(dx)=\frac 12(\delta_{-1}(dx)+\delta_{+1}(dx)).$$
This yields
$$p_r^{(-1)}(x,y)=H_0^{(-1)}(x)H_0^{(-1)}(y)+
rH_1^{(-1)}(x)H_1^{(-1)}(y)=
1+rxy,$$
giving as transition probabilities
$$\ck_{s,t}^{(-1)}(x,dy)=\frac 12(1+\frac {c(s,t)}{c(s,s)c(t,t)}xy)
(\delta_{-\sqrt{c(t,t)}}(dy)+\delta_{+\sqrt{c(t,t)}}(dy)).$$
1) fermionic Brownian motion: $X_t$ can only assume the values
$+\sqrt t$ and $-\sqrt t$ and the transition probabilities are
given by the table
$$\matrix
\kst&\sqrt t& -\sqrt t\\
\sqrt s & \frac 12(1+\sqrt{s/t})&\frac 12(1-\sqrt{s/t})\\
-\sqrt s & \frac 12(1-\sqrt{s/t})&\frac 12(1+\sqrt{s/t})
\endmatrix.$$
This case coincides with the corresponding $c=-1$ case of the
Az\'ema martingale, see \cite{Par1}.
\newline
2) fermionic Ornstein-Uhlenbeck process: This stationary process
lives on the two values $+1$ and $-1$ with the following transition
probabilities
$$\matrix
\kst&1& -1\\
1 & \frac 12(1+e^{-(t-s)})& \frac 12(1-e^{-(t-s)})\\
-1 & \frac 12(1-e^{-(t-s)})& \frac 12(1+e^{-(t-s)})
\endmatrix.$$
This classical two state Markov realization of the corresponding
fermionic relations has been known for a long time, see \cite{FB}.
\newline
3) fermionic Brownian bridge:
$X_t$ can only assume the values
$+\sqrt{t(1-t)}$ and $-\sqrt{t(1-t)}$ and 
the transition probabilities are
given by the table
$$\matrix
\kst&\sqrt{t(1-t)} & -\sqrt{t(1-t)}\\
\sqrt{s(1-s)}
 & \frac 12(1+\sqrt{\frac{s(1-t)}{t(1-s)}})&
\frac 12(1-\sqrt{\frac{s(1-t)}{t(1-s)}})\\
-\sqrt{s(1-s)}
 & \frac 12(1-\sqrt{\frac{s(1-t)}{t(1-s)}})&
\frac 12(1+\sqrt{\frac{s(1-t)}{t(1-s)}})
\endmatrix.$$
\endexample
\example{4.11. Example: Hypercontractivity}
Consider the $q$-Ornstein-Uhlenbeck process with stationary
transition operators 
$\cKqt:=\cK^{qOU}_{s,s+t}$.
Note that this $q$-Ornstein-Uhlenbeck semigroup is nothing
but the second quantization of the simplest contraction,
namely with the one-dimensional real Hilbert space $\HR=\RR$ and
the corresponding identity operator $\id:\RR\to\RR$ we
have
$$\Gamma_q(\RR)\cong L^\infty(-2/\sqrt{1-q},2/\sqrt{1-q},
\nu_q(dx))\qquad\text{and} \qquad
\Gamma_q(e^{-t}\id)\cong \cKqt.$$
We have seen that the $\cKqt$ are,
for all $t>0$, contractions on $L^2$ and on $L^\infty$ (and thus,
by 
duality and interpolation, on all $L^p$). In the classical case $q=1$ 
(and also for $q=-1$) it
is known \cite{Sim1,Nel1,Nel2,Gro,CL}
that much more is true, namely the Ornstein-Uhlenbeck
semigroup is also hypercontractive, i.e. it is bounded as a map
from $L^2$ to $L^4$ for sufficiently large $t$. Having the concrete
form of the kernel 
$$\ckqt(x,dy)=p^{(q)}_{e^{-t}}(x,y)\nu_q(dy)$$
of $\cKqt$ it is easy to check that
we also have hypercontractivity for all $-1<q<1$. Even more,
we can show that $\cKqt$ is bounded from $L^2$ to $L^\infty$ for
$t>0$, i.e. we have what one might call
\lq ultraconctractivity' -- which is,
of course, not given for $q=\pm 1$. This ultracontractivity 
follows from the estimate
$$\Vert \cKqt h\Vert_\infty\leq \alpha(t,q)^{1/2} \Vert h\Vert_2
\qquad\text{where}\qquad
\alpha(t,q):=\sup_{x\in\lb -2,2\rb}\sup_{y\in\lb -2,2\rb} 
p_{e^{-t}}^{(q)}(x,y)$$
and from the explicit form of $p_r^{(q)}$ from Theorem 1.10,
which ensures that $\alpha(t,q)$ is finite for $t>0$ and
$-1<q<1$.
One may also note that for small $t$ the leading term 
of $\alpha(t,q)^{1/2}$ is of order $t^{-3/2}$.
\endexample
\example{4.12. Open Problems}
1) The situation concerning classical versions of non-Markovian
$q$-Gaussian processes is not clear at the moment.
\newline
2) Consider a symmetric measure $\mu$ on $\RR$ with compact support.
Then there exist a sequence of polynomials $(P_n)_{n\in\Nn}$ such
that $P_n$ is of degree $n$ and such that these polynomials are 
orthogonal with respect to $\mu$. Let us define a semigroup
$U_t$ on $L^2(\mu)$ by 
$$U_t P_n=e^{-nt}P_n.$$
If these $U_t$ are positivity
preserving then they constitute the transition
operators of a stationary Markov process,
whose stationary distribution is given 
by $\mu$. Our $q$-Ornstein-Uhlenbeck
process is an example of this general construction for the
measure $\nu_q$. 
The existence of the functor $\sqq$ \lq explains' the fact that
the $q$-Ornstein-Uhlenbeck semigroup 
is positivity preserving from a more general
(non-commutative) point of
view --
note that although Theorem 2.11 is for $\dim\HR=\dim \HR'=1$ a
purely commutative statement, its proof is even in this case
definitely non-commutative.
Of course, not for all measures $\mu$ the semigroup
$U_t$ is positivity preserving.
But one might wonder whether it is possible to find
for each measure with this property 
--or at least for some special class of such measures --
some analogous kind of
functor. See also \cite{BSp4} for related investigations.
\endexample

\enddocument